\documentclass[letterpaper, 10 pt, conference]{ieeeconf}  

\IEEEoverridecommandlockouts    
\usepackage[USenglish]{babel}

\usepackage{xcolor,graphics,graphicx} 
\usepackage{epsfig} 
\usepackage{amsmath} 
\usepackage{amssymb}  
\usepackage{amsfonts}
\usepackage{mathtools}  
\usepackage{optidef,booktabs}
\usepackage{cite}
\usepackage{url}
\usepackage{lipsum} 
\usepackage{hyperref}
\usepackage[utf8]{inputenc} 
\usepackage[font=small,skip=0pt]{subcaption}
\usepackage[font=small,skip=0pt]{caption}
\usepackage{tabularray}
\usepackage{tikz}
\usetikzlibrary{decorations.markings}

\providecommand{\SRG}[1]{\operatorname{SRG}(#1)}

\UseTblrLibrary{amsmath}

\newtheorem{remark}{Remark}

\newtheorem{theorem}{Theorem}
\newtheorem{corollary}{Corollary}

\title{\LARGE \bf Mixed Small Gain and Phase Theorem: \\A new view using Scale Relative Graphs}

\author{Eder~Baron-Prada, Adolfo Anta, Alberto Padoan and~Florian Dörfler
\thanks{Eder Baron is with the Austrian Institute of Technology, 1210 Vienna, Austria, and also with the Automatic Control Laboratory, ETH Zurich, 8092 Z\"urich, Switzerland. (e-mail: ebaron@ethz.ch)}
\thanks{Adolfo Anta is with the Austrian Institute of Technology, Vienna 1210, Austria (e-mail: adolfo.anta@ait.ac.at).}
\thanks{Alberto Padoan is with the Department of Electrical \& Computer Engineering, University of British Columbia. (e-mail: apadoan@ece.ubc.ca).}
\thanks{Florian Dörfler are with the Automatic Control Laboratory, ETH Zürich, Zürich 8092, Switzerland (e-mail: apadoan@ethz.ch, dorfler@ethz.ch).}
\thanks{We thank Verena Häberle, Linbin Huang, Xiuqiang He, the IfA group for the fruitful discussions, and the reviewers for their feedback and suggestions.} 
}

\begin{document}

\maketitle
\thispagestyle{empty}
\pagestyle{empty}

\begin{abstract}
We introduce a novel approach to feedback stability analysis for linear time-invariant (LTI) systems, overcoming the limitations of the sectoriality assumption in the small phase theorem. While phase analysis for single-input single-output (SISO) systems is well-established, multi-input multi-output (MIMO) systems lack a comprehensive phase analysis until recent advances introduced with the small-phase theorem.

A limitation of the small-phase theorem is the sectorial condition, which states that an operator's eigenvalues must lie within a specified angle sector of the complex plane. We propose a framework based on Scaled Relative Graphs (SRGs) to remove this assumption. We derive two main results: a graphical set-based stability condition using SRGs and a small-phase theorem with no sectorial assumption. These results broaden the scope of phase analysis and feedback stability for MIMO systems.
\end{abstract}

\section{Introduction}

In classical frequency domain analysis of SISO LTI systems, the magnitude response and phase response play important roles\cite{zhou1998}. Magnitude and phase can be visualized in Bode and Nyquist diagrams, representing the system behavior in response to inputs across different frequencies\cite{1102280}. This well-established theory covers SISO systems comprehensively, while the analysis becomes more complex when extended to MIMO systems\cite{skogestad2005}. For MIMO systems, the focus has historically been placed on the gain, supported by the development of robust control strategies and the widely used small gain theorem\cite{zhou1998}. However, the phase component has been less developed in MIMO cases. Early approaches to address this gap include the principal region approach \cite{Postlethwaite1981} approach and the mixed property approach\cite{GRIGGS_2012_sufficientNyquist}.

Recognizing this gap, recent efforts, such as the work by Chen et al. \cite{Chen2020}, have aimed to extend the theory of phases for MIMO systems. They propose a framework\cite{Wang2024}, that advances the concept of complex operator phases\cite{Chen2022}. The second part of this trilogy stands out by offering an extensive phase theory for MIMO systems~\cite{Chen2019}. This includes the small phase theorem, a counterpart of the small gain theorem, used in feedback stability analysis.

The small phase theorem marks a significant advancement for LTI systems. Rooted in the concept of numerical range, first applied in control theory by Owens \cite{Owens1984}, this theorem complements the existing small gain theorem, offering a new lens for analyzing stability. However, applying the small phase theorem requires a condition called sectoriality \cite{Wang2024,Chen2019,Chen2020}, which is often not fully satisfied in practice (e.g. power converters\cite{huang2024gain,Woolcock2023}).

To address the limitations imposed by the sectoriality assumption, we propose a novel approach based on the concept of SRG, initially introduced in optimization theory \cite{ryu2022large,Ryu_2021}, and later used for stability analysis of nonlinear operators \cite{Chaffey2021,Chaffey_2023,CHAFFEY_2022}. Using SRGs, we derive a flexible stability condition that avoids the sectoriality constraint, making it applicable to a wider range of systems.

Our contributions are twofold: first, we use the SRG definition and the principle of superposition to develop a graphical, set-based stability criterion. The second result presents an improved small-phase condition that eliminates the sectorial condition, offering a less conservative alternative to the small-phase theorem. We use this theorem to prove a new version of the mixed small gain and phase theorem.


\section{Preliminaries}\label{sec:basics}

The sets of real and complex numbers are denoted by $\mathbb{R}$ and $\mathbb{C}$, respectively. The polar representation of $z\in \mathbb{C}$ is defined as $z=re^{j\alpha}$ with $r$ denoting the gain and $\alpha$ the angle. When referring to the angle of $z$, we use $\angle(z)$. The complex conjugate of $z \in \mathbb{C}$ is denoted by $z^*$, and its real part is denoted as $\Re(z)$. A set $S\subseteq \mathbb{C}$ is said to be convex if for all $0 \leq t \leq 1$ and $s_1, s_2 \in S$, then  $t s_1 + (1-t) s_2 \in S$. We denote the imaginary unit as $j$. The time derivative of $x$ is denoted as $\dot{x}$. An operator $ A $ is invertible if there exists an operator $ A^{-1} $ such that $ AA^{-1} = A^{-1}A = I $, where $ I $ is the identity operator. Let $\mathcal{H}$ denote a Hilbert space defined over the field $F$. An operator $A: \mathcal{H} \rightarrow \mathcal{H}$ is linear if $A(\alpha x + \beta y) = \alpha A(x) + \beta A(y)$ such that for any $\alpha, \beta \in F$ and $x, y \in \mathcal{H}$. The spectrum of $ A $ consists of all scalar values $ \lambda_i \in \mathbb{C}$ such that $ (A - \lambda_i I) $ is not invertible.  

\subsection{Signal Spaces}

We focus on Lebesgue spaces of square-integrable functions, $\mathcal{L}_2$. Given the time axis, $\mathbb{R}_{\ge 0}$, and a field $F \in \{\mathbb{R}, \mathbb{C}\}$, we define the space $\mathcal{L}_2^n(F)$ by the set of signals $u: \mathbb{R}_{\ge 0} \rightarrow F^n$ and $y: \mathbb{R}_{\ge 0} \rightarrow {F}^n$ such that the inner product of $u, y \in \mathcal{L}_2^n({F})$ is defined by $\langle u, y \rangle := \int_0^{\infty} u(t)^* y(t) \, dt$, and the norm of $u$ is defined by $\|u\| := \sqrt{\langle u, u\rangle }$. The Fourier transform of $u \in \mathcal{L}_2^n({F})$ is defined as $\hat{u}(j\omega) := \int_0^{\infty} e^{-j\omega t} u(t) \, dt$.

\subsection{Transfer Functions and Stability of LTI Systems} 

Transfer functions describe the input-output behavior of LTI systems, defined by 
\begin{align*}
    \dot{x} = Ax + Bu;\qquad y = Cx + Du
\end{align*} where $x \in \mathbb{R}^n$ is the state vector, $u \in \mathbb{R}^m$ is the input, and $y \in \mathbb{R}^p$ is the output, with appropriately dimensioned matrices $A$, $B$, $C$, and $D$. This work considers the space $\mathcal{RH}_\infty$ of rational, proper, stable transfer functions representing bounded, causal LTI operators. Such systems induce an input-output gain that quantifies the relative output size to the input. For $\mathcal{L}_2$ signals, this gain corresponds to the $H_\infty$ norm\cite{Chen2022,zhou1998,arcozzi2020hardy}. An LTI system is $\mathcal{L}_2$-stable if a bounded input $u(t) \in \mathcal{L}_2$ produces a bounded output $y(t) \in \mathcal{L}_2$.
 \subsection{Review: Small phase theorem}
 In this subsection, we introduce the small phase theorem following \cite{Chen2020,Chen2019,Chen2022,Wang2024} which provides a framework for phase-based stability analysis in LTI systems. The numerical range of an operator $ A \in \mathbb{C}^{n\times n}$ is defined by \cite{psarrakos2002numerical}
 \begin{align*}
     W(A) = \{ \langle Ax, x \rangle : x \in \mathbb{C}^n, \|x\| = 1 \}. 
 \end{align*}
 
$ W(A) $ is a convex subset of $ \mathbb{C} $ and contains the eigenvalues of $ A $. An operator is called \textit{sectorial} if $ W(A) $ is contained within an angular sector of the complex plane, and $ 0 \notin W(A) $ as shown in Fig.~\ref{fig:Sectorial}. The \textit{sectorial decomposition} of a sectorial operator $ A $ is defined as $A = T D T^{-1}$, where $ T $ is an invertible operator and  $ D $ is a unique diagonal unitary operator\cite{Chen2019,Wang2024}. In addition, the elements of $ D $ are the eigenvalues of $ A $, and lie on an arc of the unit circle of length smaller than $ \pi $. The phases of the operator $ A $ are defined as the angle of each element of $ D $, and denoted as
\begin{align*}
\alpha_{\max}(A)= \alpha_1(A)\geq \ldots\geq\alpha_n(A) = \alpha_{\min}(A).
\end{align*}

The phases are contained in $ \alpha_{\max}(A)-\alpha_{\min}(A)<\pi$. Besides, the supporting angles are defined as the maximum and minimum eigenvalues angles\cite{Wang2024}. Note that the sectorial decomposition is necessary for the computation of the operator phases\footnote{For a more detailed explanation of the phase computation of a sectorial operator, we direct interested readers to \cite{Wang2024,Chen2019}.}.

\subsubsection{Classification of Sectorial Operators}

\begin{figure}[h!]
    \centering
\begin{subfigure}[b]{0.15\textwidth}
    \centering
    \includegraphics[width=1\linewidth]{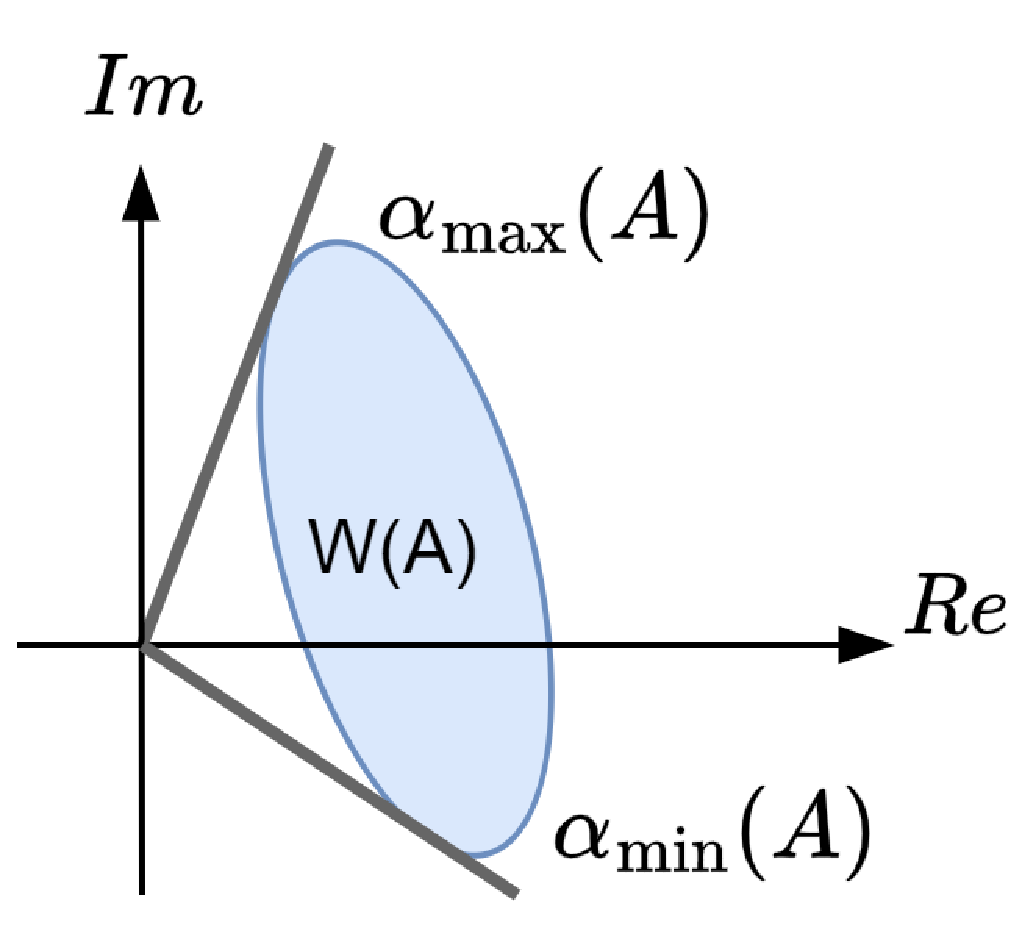}
    \caption{} \label{fig:Sectorial}
\end{subfigure}    
\begin{subfigure}[b]{0.15\textwidth}
    \centering
    \includegraphics[width=1\linewidth]{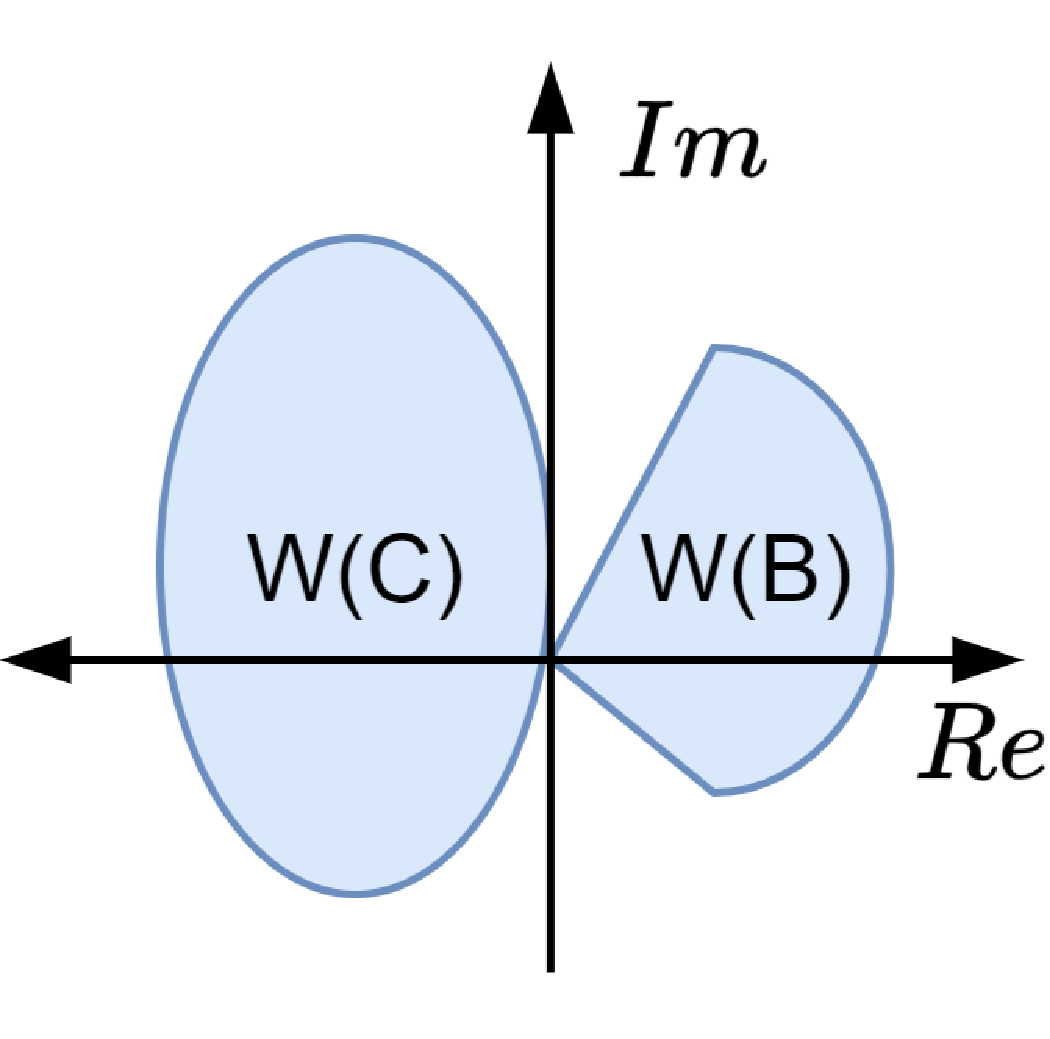}
    \caption{}     \label{fig:semi-quasi}
\end{subfigure}
\begin{subfigure}[b]{0.15\textwidth}
    \centering
    \includegraphics[width=1\linewidth]{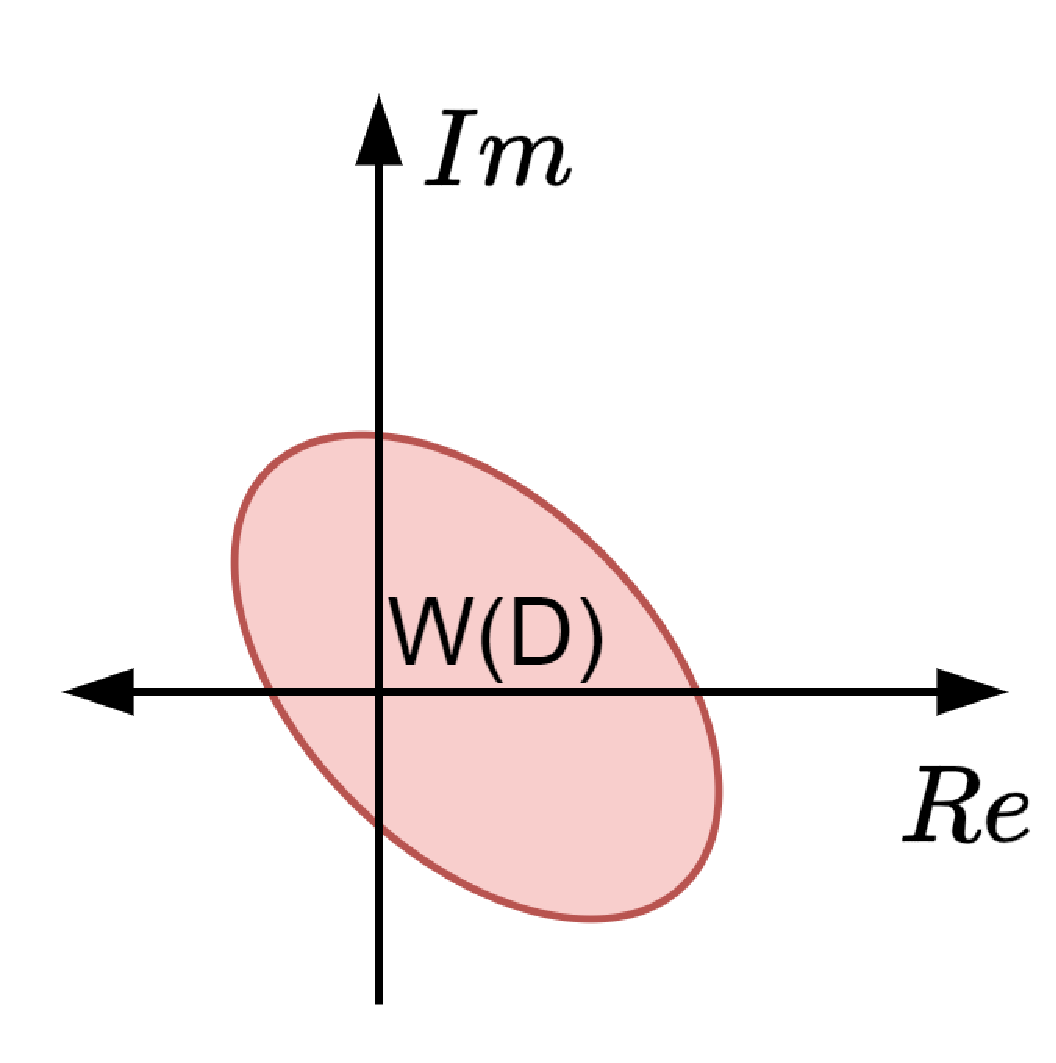}
    \caption{} \label{fig:nosectorial}
\end{subfigure}    
    \caption{(a) $W(A)$ of a sectorial operator $A$ with supporting angles ${\alpha}_{\max}(A)$ and ${\alpha}_{\min}(A)$. (b) $W(B)$ of a quasi-sectorial operator $B$ and $W(C)$ of a semi-sectorial operator $C$, respectively. (c) $W(D)$ of a non-sectorial operator $D$.}
    \label{fig:sectorialfigures}
\end{figure}

Operators are classified based on the properties of their numerical range\cite{Chen2020}. An operator is \textit{quasi-sectorial} if the supporting lines of $ W(A) $ form an angle smaller than $ \pi $ and $ 0 \in W(A) $. An operator is \textit{semi-sectorial} if the supporting lines of $ W(A) $ form an angle less than or equal to $ \pi $, and $ 0 \in W(A) $ as shown in Fig.\ref{fig:semi-quasi}. Finally, an operator is non-sectorial if its numerical range includes 0 in its interior as in Fig.\ref{fig:nosectorial}.
 
\subsubsection{The Small Phase Theorem }
The small phase theorem proposed in \cite{Wang2024,Chen2019,Chen2022,Chen2020}, provides a framework for extending phase analysis from SISO to MIMO systems, while also broadening the scope of the passivity theorem and concepts like passivity in LTI systems\cite{Wang2024}. However, the theorem is limited by the requirement that the system frequency response must exhibit sectoriality, a condition that real-world systems do not always fulfill\cite{huang2024gain,Woolcock2023}.
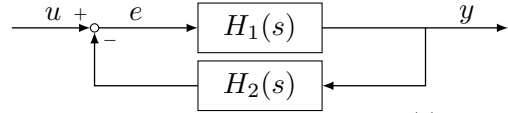
\begin{figure}[hbt]
    \centering
    \begin{tikzpicture}[scale=1.1, every node/.style={transform shape}]
\draw (9.75,4.3) rectangle (11.25,3.7);
\node at (10.5,4) {$H_1(s)$};
\draw (9.75,3.6) rectangle (11.25,3); 
\node at (10.5,3.3) {$H_2(s)$};
\draw[-latex, line width = .5 pt] (12.5,4) -- (13.5,4);
\draw[-latex, line width = .5 pt] (11.25,4) -- (12.5,4) -- (12.5,3.3) -- (11.25,3.3);
\draw[fill = white] (8.5,4) circle [radius=0.05];
\draw[-latex, line width = .5 pt] (9.75,3.3) -- (8.5,3.3) -- (8.5,3.95) ;
\draw[-latex, line width = .5 pt] (8.55,4) -- (9.75,4);
\draw[-latex, line width = .5 pt] (7.5,4) -- (8.45,4);
\node at (8.6892,3.8444) {\tiny$-$};
\node at (8.3332,4.1664) {\tiny$+$};
\node at (9,4.1664) {$e$};
\node at (8,4.1664) {$u$};
\node at (13,4.1664) {$y$};
\end{tikzpicture}
\caption{Feedback interconnection between $H_1(s)$ and $H_2(s)$. }
    \label{fig:fb}
\end{figure}
\begin{theorem}\textbf{(Small phase theorem\cite{Chen2019}) } \label{Thm:smallphase}
Assume $ H_1(s) \in \mathcal{RH}_\infty$ and $ H_{2}(s) \in \mathcal{RH}_\infty$ are connected in a feedback loop as shown in Fig. \ref{fig:fb}. If for each $s=j\omega$, with $ \omega \in [0, \infty) $, the following holds: 
\begin{enumerate}
    \item ${\alpha}_{\max}(H_1(s)) + {\alpha}_{\max}(H_{2}(s)) < \pi$, {and}
    \item ${\alpha}_{\min}(H_1(s)) + {\alpha}_{\min}(H_{2}(s)) > -\pi,$ {and} 
    \item $H_1(s)$ and $H_{2}(s)$ satisfies the sectorial property
\end{enumerate}
then, the closed-loop system is stable.
\end{theorem}
\section{ Scaled Relative Graphs (SRGs)}\label{sec:SRG}

SRGs are initially introduced in \cite{ryu2022large,Ryu_2021}. It is originally a tool used in optimization theory for convergence analysis\cite{Ryu_2021}. The SRG tool can be used not only with linear operators but also to address non-linear operators\cite{Chaffey2021,Chaffey_2023}. This section recalls a stability condition based on SRGs proposed by \cite[Theorem 1]{Chaffey_2023}. Finally, we present an SRG-based stability condition for LTI systems.
\subsection{ Scaled Relative Graphs of Square Matrices}

Consider an operator $A: \mathcal{H} \rightarrow \mathcal{H}$.  The SRG of $A$ is defined as \cite{Ryu_2021}
\begin{align*}
    \operatorname{SRG}(A)= \left\{\dfrac{\| {y_2-y_1}\|}{\| {u_2-u_1}\|} \exp \left[\pm j \angle(u_2-u_1,y_2-y_1) \right]\right\},
\end{align*}
where $ {u_1,u_2} \in \mathcal{H}$ are a pair of inputs which outputs are ${y_1,y_2}$, i.e. ${y_{1}}=A({u_{1})}$ and ${y_{2}}=A({u_{2})}$. If $y_2-y_1=0$ or $u_2-u_1=0$, $\angle(u_2-u_1,y_2-y_1)=0$\cite{Ryu_2021}. If $A$ is a square matrix, $\operatorname{SRG}(A)$ is given by  \cite{pates2021scaled}
\begin{align}
    \operatorname{SRG}(A)= \left\{\dfrac{\| {y}\|}{\| {u}\|} \exp \left( \pm j \arccos \left(\dfrac{\Re(\langle {y},  {u}\rangle)}{\| {y}\|\| {u}\|}\right)\right)\right\},
    \label{eqn:SRG_operator}
\end{align}
with $\| {u}\|=1 $, and ${y}=A  {u}$. The ratio $\frac{\|y\|}{\|u\|}$ represents the magnitude change in the output compared to the input. The term $\frac{\Re(\langle y, u \rangle)}{\|y\| \|u\|}$ measures the angle difference between the output and the input \cite{millman1993}. Furthermore, as in the nonlinear case, if $\|y\|=0$, then $\angle(u,y)=0$.  Note that, for the linear and nonlinear versions, the $\operatorname{SRG}(A)$ is symmetric with respect to the real axis. 

The arc property defines specific arcs between points $z$ and $z^*$ based on their real parts\cite{Ryu_2021}. The right-arc property, denoted as $\text{Arc}^+(z, z^*)$, is the arc between points $z$ and $z^*$ centered at the origin, where the real part of the arc is greater than or equal to $\Re(z)$, i.e., $ \text{Arc}^+(z, z^*):=\{re^{j(1-2\theta)\alpha}|z=re^{j\alpha}, \alpha \in (-\pi,\pi],\theta\in[0,1], r\geq 0\}$. Conversely, the left-arc property, $\text{Arc}^-(z, z^*)$, is similarly defined but with the real part less than $\Re(z)$, i.e., $ \text{Arc}^-(z, z^*):= -\text{Arc}^+ (-z, -z^*)$. 


\subsection{Stability Conditions based on SRGs}
Theorem \ref{thm:GFT} provides the foundation for certifying stability using SRGs in operators with finite incremental $\mathcal{L}_2$ gain.

\begin{figure}[hbt]
  \centering
    \begin{tikzpicture}[scale=1.1, every node/.style={transform shape}]
\draw (9.75,4.3) rectangle (11.25,3.7);
\node at (10.5,4) {$H(s)$};
\draw[-latex, line width = .5 pt] (12.5,4) -- (13.5,4);
\draw[-latex, line width = .5 pt] (11.25,4) -- (12.5,4) -- (12.5,3.3) -- (8.5,3.3) -- (8.5,3.95);
\draw[fill = white] (8.5,4) circle [radius=0.05];
\draw[-latex, line width = .5 pt] (8.55,4) -- (9.75,4);
\draw[-latex, line width = .5 pt] (7.5,4) -- (8.45,4);
\node at (8.6892,3.8444) {\tiny$-$};
\node at (8.3332,4.1664) {\tiny$+$};
\node at (9,4.1664) {$e$};
\node at (8,4.1664) {$u$};
\node at (13,4.1664) {$y$};
\end{tikzpicture}
    \caption{Unitary feedback connection of  $H(s)$. }
    \label{fig:unitary_fb}
\end{figure}
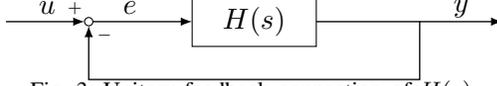

\begin{theorem}\textbf{\cite[Theorem 1]{Chaffey_2023}}\label{thm:GFT}
Assume $ H : \mathcal{L}_2 \rightarrow \mathcal{L}_2$ is an operator with finite incremental $\mathcal{L}_2$ gain.  If 
\begin{align}
    (-1,0) \notin \tau\operatorname{SRG}(H),\;\; \forall\tau \in(0,1].
    \label{eqn:SRGequi}
\end{align}

Then, the closed loop operator shown in Fig.\ref{fig:unitary_fb} is $\mathcal{L}_2 $ stable.
\end{theorem}
It is important to note that Theorem \ref{thm:GFT} guarantees only $\mathcal{L}_2$ stability and does not ensure finite incremental gain. If finite incremental gain is required, one needs the distance between $\SRG{H}$ and $-1$ to be non-zero $\forall \tau$\cite{chaffey2024homotopy}.
In what follows, we use Theorem \ref{thm:GFT} to derive an SRG-based stability condition for LTI systems. This is a direct corollary of Theorem \ref{thm:GFT} using linearity and time invariance.
\begin{corollary}\label{corollary:nyquist}  Assume $H(s) \in \mathcal{RH}_\infty$, if for each $s=j\omega$ with $\omega \in [0, \infty)$ 
\begin{align}
    (-1,0) \notin \tau\operatorname{SRG}(H(s)),\;\; \forall \tau \in (0,1]\label{eqn:proposition_nyquist}
    \end{align}
     Then, the closed-loop system is $\mathcal{L}_2$ stable.
\end{corollary} 
%

Corollary~\ref{corollary:nyquist} shows that when Theorem \ref{thm:GFT} is applied to LTI systems, stability can be assessed by comparing the system SRG at each frequency $\omega \in [0, \infty)$ against the point $(-1,0)$. The key advantage of Corollary \ref{corollary:nyquist} lies in leveraging the superposition principle\cite{zhou1998}. This principle allows us to treat the system response as a collection of operators acting independently at each frequency. In general, Corollary \ref{corollary:nyquist} is more conservative than the Generalized Nyquist Criterion. This is because Corollary \ref{corollary:nyquist} imposes a condition that is evaluated frequency by frequency, whereas the Generalized Nyquist Criterion incorporates two distinct criteria: one that is frequency-wise and another that accounts for the trajectory of the eigenvalues over the entire frequency spectrum \cite{Baron2025SRG}.

\section{Mixed Gain and Phase Theorem}\label{sec:mixed}

We derive stability conditions using SRGs, following the principles of the small gain and small phase theorems \cite{Chen2020}. While SRGs provide a rigorous framework, their practical use is limited by significant computational challenges. SRGs are closed-bounded sets requiring numerous points for accurate representation at each frequency, making them computationally intensive, especially for high-dimensional or real-time systems. In contrast, this section presents sufficient conditions simplifying this burden by focusing on two key characteristics at each frequency: $\sigma_{\max}(\cdot)$, related to the maximum gain, and $\alpha_{\max}(\cdot)$, related to the maximum angle between the input and the output, also called phase. This approach only requires the boundary of the SRG to calculate such points. Hence, it offers a more manageable and efficient way to analyze system stability, but at the expense of being a more conservative stability condition.
 
\subsection{Maximum gain and phase} \label{subsec:SVD_SAD}

We recall the SRG definition in \eqref{eqn:SRG_operator}, which has two main components: the ratio between input and output norms indicating magnitude change and the angle between input and output. We identify the maximum values of these components: the highest gain, $\sigma_{\max}(\cdot)$, and the maximum phase, $\alpha_{\max}(\cdot)$.

\subsubsection{Maximum Gain}
the gain calculation between the input and output vector is given by:
\begin{align}
\sigma_{\max}(A)=\max_{||u||=1} \frac{||Au||}{||u||} \label{eqn:SingularValuedef}.
\end{align}

Note that \eqref{eqn:SingularValuedef} defines the maximum singular value, a key component of the small gain theorem \cite{zhou1998}.

\subsubsection{Maximum phase}
the angle calculation between the input and output vector is given by:
\begin{align}
\hat{\alpha}_{\max}(A) :=\max_{||Au||\neq0, \|u\|=1} \arccos \left(  \dfrac{\Re(\langle Au,u \rangle)}{||Au||||u||} \right).\label{eqn:SingularAngledef}
\end{align}
If $\|Au\|=0\;\forall\; \{u|\; \|u\|=1\}$, then $\hat{\alpha}_{\max}(A)~=~0$. In addition, note that the maximum gain, $\sigma_{\max}(\cdot)$, and maximum phase, $\hat{\alpha}_{\max}(\cdot)$, typically do not occur at the same point on the SRG for MIMO systems, whereas for SISO systems, it is always the same point.

 \subsection{Small phase theorem using SRGs}
We now offer an alternative to Theorem \ref{Thm:smallphase}.  Our approach does not require sectorial properties, dropping the strongest constraint in Theorem \ref{Thm:smallphase}. 
\begin{theorem}
\textbf{(Small phase theorem based on SRGs)} \label{thm:Small_phase_SRG}
Assume $ H_1(s) \in \mathcal{RH}_\infty $ and $ H_2(s) \in \mathcal{RH}_\infty$. If for each $s=j\omega$, with $ \omega \in [0, \infty) $ and
\begin{align}
    \hat{\alpha}_{\max}(H_1(s)) + \hat{\alpha}_{\max}(H_2(s)) < \pi, 
    \label{eqn:SRGsmallphase}
\end{align}
Then, the closed-loop system is $\mathcal{L}_2$ stable.
\end{theorem}
\proof Proof in the Appendix \ref{proof:Small_phase_SRG}.

A similar version of Theorem \ref{thm:Small_phase_SRG} appears in \cite{chen2021singularanglenonlinearsystems}, based on the singular angle concept and computed via the matrix normalized numerical range. Even though we require the right-arc property approximation in the proof, we do not include it as an assumption in Theorem \ref{thm:Small_phase_SRG}, since over-approximating any SRG to have the right-arc property can be done without introducing any conservatism, i.e., without affecting $\alpha_{\max}(A)$, as shown later in Remark \ref{remark:Overapproximation}.
Note that Theorem~\ref{thm:Small_phase_SRG} requires only condition \eqref{eqn:SRGsmallphase}, in contrast to Theorem~\ref{Thm:smallphase}. Specifically, we can omit condition 2), which states ${\alpha}_{\min}(H_1(s))~+~{\alpha}_{\min}(H_2(s))~>~-~\pi$, because the SRG is symmetric. This symmetry makes condition 2) equivalent to ${-\hat{\alpha}}_{\max}(H_1(s)) - {\hat{\alpha}}_{\max}(H_2(s)) > -\pi$, which is equivalent to \eqref{eqn:SRGsmallphase}. Finally, condition 3), the sectorial condition, is omitted because our approach does not require this property for phase calculations. 

\begin{remark}[\textbf{Over-approximation via right-arc property}]
\label{remark:Overapproximation}
Any operator may be over-approximated by an operator with the right-arc property. It is possible to find an approximation that does not modify $\sigma_{\max}(\cdot)$ and $\hat{\alpha}_{\max}(\cdot)$. Consider an operator $A$ and denote $\operatorname{SRG}(\Tilde{A})$, as the over-approximation of $ \operatorname{SRG}(A)$ that has the right-arc property. This approximation can be found by taking each point $z\in \operatorname{SRG}(A)$ and including into $\operatorname{SRG}(\Tilde{A})$ every point in the arc defined by the right-arc property, i.e., 
\begin{align}
     \operatorname{SRG}(\Tilde{A}):=
     \{r&e^{j(1-2\theta)\alpha}|\forall z=re^{j\alpha}\in \operatorname{SRG}(A),\nonumber\\& \alpha \in [-\pi,\pi],\theta\in[0,1],\infty > r> 0\}. \label{eqn:SRG_rightarc}
\end{align}
Note that the resulting $\operatorname{SRG}(\Tilde{A})$ includes more points in the complex plane if $\operatorname{SRG}(A)$ does not have the right-arc property. In other words, $\operatorname{SRG}(\Tilde{A})\supseteq \operatorname{SRG}(A)$. Finally, it is possible to conclude from \eqref{eqn:SRG_rightarc} that $\sigma_{\max}(A)$ and $\hat{\alpha}_{\max}(A)$ remain unchanged.
\end{remark}

The small gain theorem provides a stability condition for feedback systems by ensuring the product of maximum gain is less than one. We recall it in Theorem~ \ref{theorem:smallgain}\cite{zhou1998}.
\begin{theorem}\textbf{(Small gain theorem)\cite{zhou1998}} \label{theorem:smallgain} Assume $ H_1(s) \in \mathcal{RH}_\infty $ and $ H_2(s) \in \mathcal{RH}_\infty $. If for each $s=j\omega$, with $ \omega \in [0, \infty) $ and
\begin{align}
      \sigma_{\max}(H_1(s))\sigma_{\max}(H_2(s)) < 1, 
     \label{eqn:SRGsmallgain}
\end{align}
Then, the closed-loop system is $\mathcal{L}_2$ stable.
\end{theorem}

We can propose the following mixed gain and phase theorem. Even though this theorem has been previously established\cite{Chen2019,huang2024gain,Woolcock2023}, in our new framework, we can restate it as follows. 

\begin{theorem}\textbf{(Mixed phase and gain Theorem)}\label{thm:mixedpg}\\
Assume $H_{1}(s) \in \mathcal{RH}_\infty$ and $H_{2}(s)\in \mathcal{RH}_\infty$. If for each $s=j\omega$ with $\omega \in [0,\infty)$, {either}
\begin{enumerate}
    \item The \textbf{phase condition} holds, i.e., 
\begin{align*}
    \hat{\alpha}_{\max}(H_{1}(s))+\hat{\alpha}_{\max}(H_{2}(s))< \pi, \textbf{or}
\end{align*}
\item  The \textbf{gain condition} holds, i.e.,
\begin{align*}
\sigma_{\max}(H_{1}(s)) \sigma_{\max}(H_{2}(s))<1 
\end{align*}
\end{enumerate}
Then, the closed-loop system shown in Fig. \ref{fig:fb} is stable.
\end{theorem}

\proof  Proof in the Appendix \ref{proof:mixedpg}

If the phase condition holds across the entire spectrum, we can conclude that a system is passive \cite{Chaffey_2023,Wang2024}. Additionally, if the gain condition is satisfied throughout the frequency spectrum, the system is contractive \cite{baron2023parameter}. 

\subsection{Comparison between different phase calculations}

The definition of phase varies between the numerical range-based and the SRG-based phase. In our approach, we utilize the SRG to represent, at each frequency, the set of potential phases that a system can exhibit. These phases correspond to the angles of the image of the unit sphere under the linear transformation of a square matrix $A\in\mathbb{C}^{n\times n}$. Specifically, they are determined by the angular difference of $y$ and $x$, where $y= Ax$ , i.e., $\arccos \left( \frac{\Re(\langle y, x \rangle)}{\|y\| \|x\|} \right),$ 
where $x \in \mathbb{C}^n$  and $\|x\| = 1$. By contrast, the phase definition used in \cite{Wang2024,Chen2019} is derived from the image of the unit sphere of the quadratic form of $A$, where the phases are obtained as the angular component of $y = x^\top A x$ for all $x \in \mathbb{C}^n $ with $\|x\| = 1$. Consequently, the phase values calculated in these two approaches generally do not coincide.

\section{Numerical Examples}\label{sec:examples}

\subsection{Comparison between different Small-Phase Theorems}
\begin{figure}[h]
    \centering
    \includegraphics[width=1\linewidth]{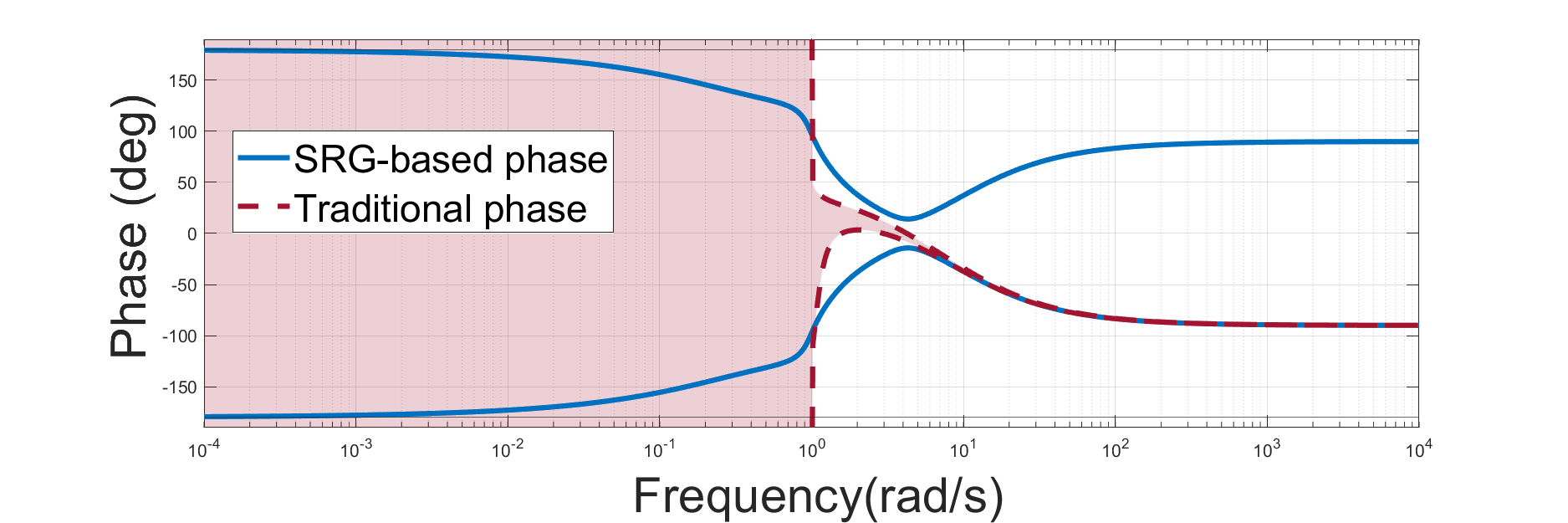}
    \caption{Comparison between SRG-based phase calculation and phase calculation using sectorial properties.}
    \label{fig:comparison}
\end{figure}

\begin{figure}[h]
    \centering 
\begin{subfigure}[b]{0.235\textwidth}
    \centering
    \includegraphics[width=0.8\linewidth]{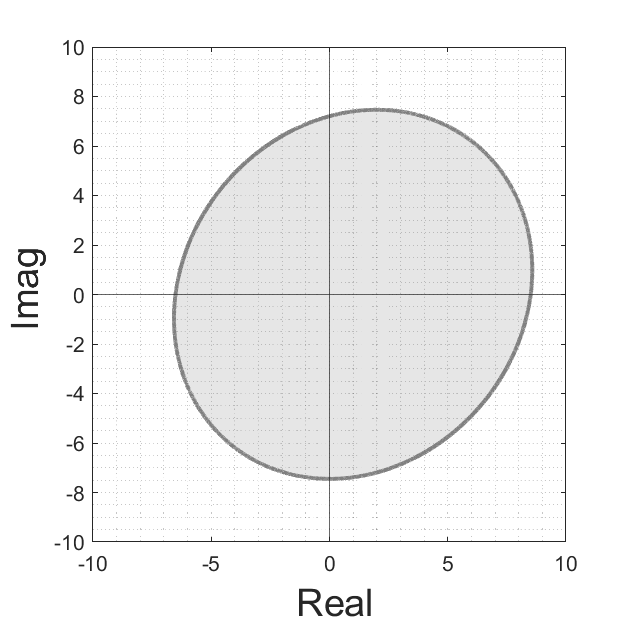}
    \caption{}
\end{subfigure}    
\begin{subfigure}[b]{0.235\textwidth}
    \centering
    \includegraphics[width=0.8\linewidth]{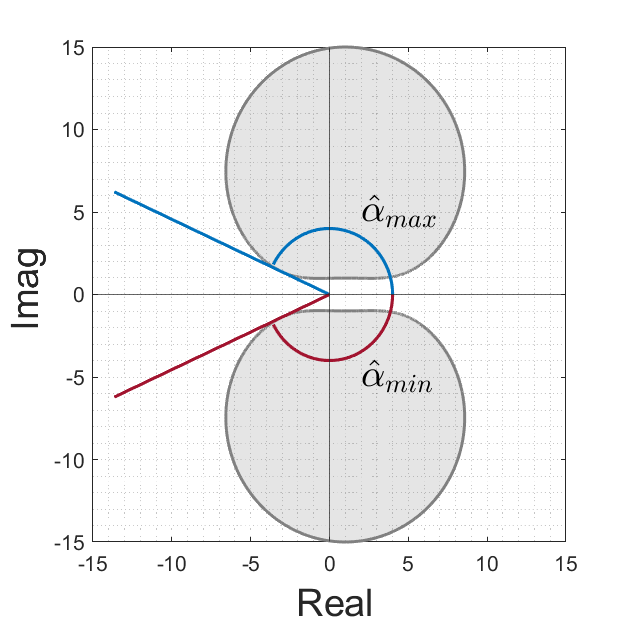}
    \caption{}     
\end{subfigure}  
    \caption{(a) Numerical range of $H_1(j\omega)$ in gray for $\omega=0.1$ rad/s (b) $\SRG{H_1(j\omega)}$, with $\hat{\alpha}_{\max} (H_1(j\omega))$ in blue and $\hat{\alpha}_{\min}(H_{1}(j\omega))$ in red with $\omega=0.1$ rad/s. }
    \label{fig:comparison_one}
\end{figure}

Consider $H_1(s)$ and $H_2(s)$ as 
\begin{align*}
H_1(s)&=\left[\begin{array}{cc} \frac{20s+30}{s^2+13\,s+30} & \frac{10}{s^2+11s+10}\\ \frac{-15}{s^2+10s} & \frac{20s^2+40\,s+30}{s^3+14s^2+43s+30}  \end{array}\right],\\
H_2(s)&=\left[\begin{array}{cc} \frac{50s+2500}{s^2+100s+2501} & 
\frac{-50}{s^2+100s+2501}\\ \frac{50}{s^2+100s+2501} & \frac{50s+2500}{s^2+100s+2501}\end{array}\right],
\end{align*}
in feedback as in Fig.~\ref{fig:fb}. Initially, we compare the phases calculated for $H_1(s)$ by each approach shown in Fig. \ref{fig:comparison}. In blue, the maximum and minimum SRG-based phases are depicted, which are symmetric, i.e., $\hat{\alpha}_{\max}(H_{1}(j\omega))=-\hat{\alpha}_{\min}(H_{1}(j\omega))$ across the entire frequency spectrum. The numerical range-based maximum and minimum phases are shown as dashed red line, calculated as in \cite{Chen2019,Chen2020}. Note that for frequencies below $\omega<1$  rad/s, the transfer function $H_1(s)$ is not sectorial as can be seen in Fig.\ref{fig:comparison_one} for $\omega=0.1$ rad/s. It is possible to see in Fig.\ref{fig:comparison} that from $\omega>10$ rad/s, the minimum phases are the same for both approaches.
  
\begin{figure}[ht]
    \centering 
    \includegraphics[width=1\linewidth]{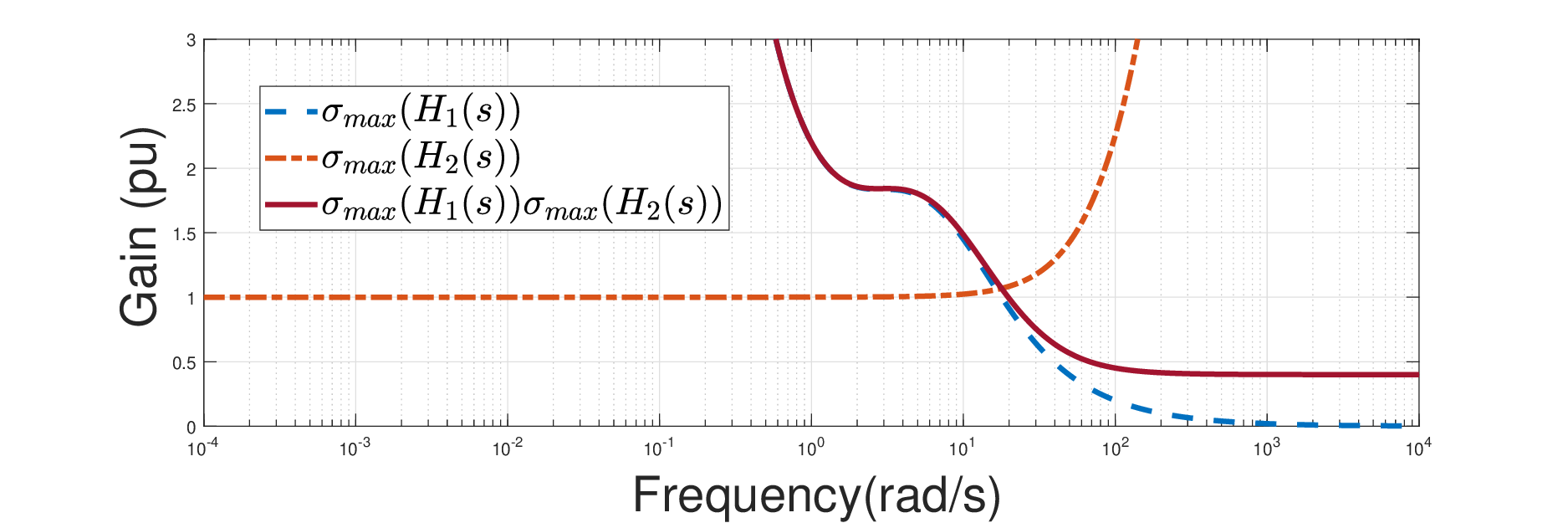}
    \caption{Gain plot}     \label{fig:smallgain}
    \includegraphics[width=1\linewidth]{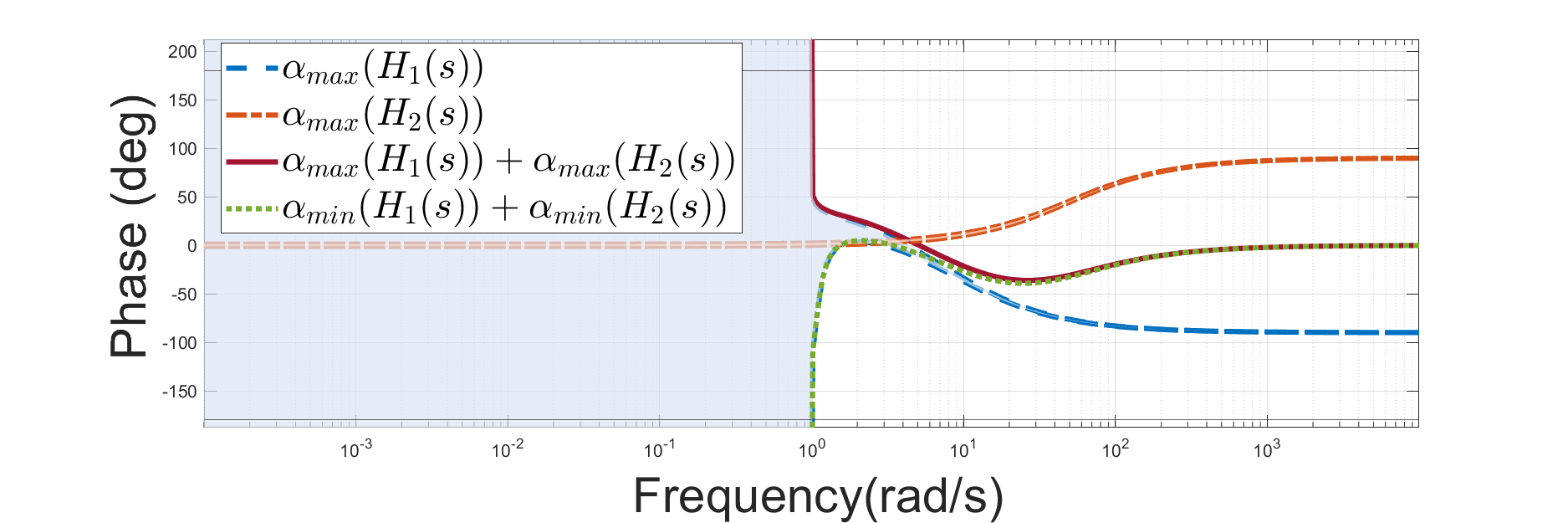}
    \caption{Traditional phase plot}         \label{fig:smallphase_traditional}
    \includegraphics[width=1\linewidth]{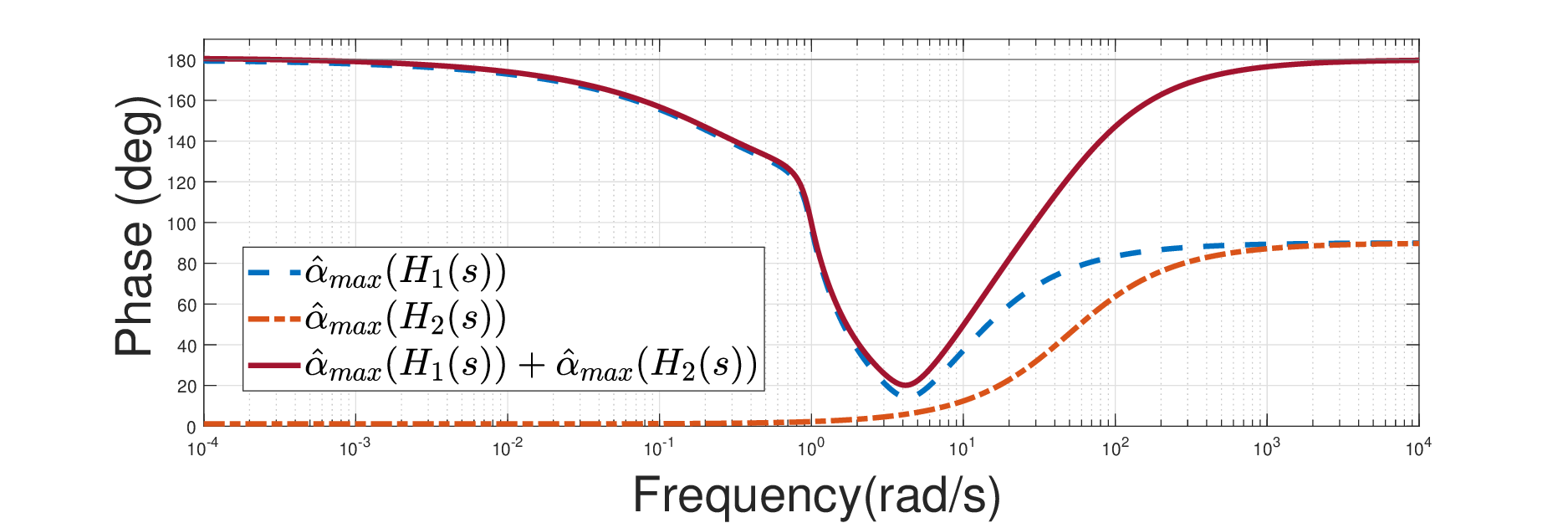}
    \caption{SRG-based phase plot }     \label{fig:smallphaseSRG}
\end{figure}

We use Theorem \ref{thm:mixedpg} to evaluate the feedback loop stability. Fig. \ref{fig:smallgain} shows the gain plot of $H_1(s)$ and $H_2(s)$, which reveals that Theorem~\ref{theorem:smallgain} certifies the system stability for $\omega>19.74$~  rad/s. Fig.\ref{fig:smallphase_traditional} shows the phase of both systems calculated as \cite{Chen2020}. The system is stable for $\omega>1$  rad/s. For $\omega<1$ rad/s, $H_1(s)$ is not sectorial; therefore, no conclusion can be reached.
 
On the other hand, the SRG-based phase is shown in Fig.\ref{fig:smallphaseSRG}. Note that as our interest is to evaluate Theorem \ref{thm:Small_phase_SRG}, it suffices to plot $\hat{\alpha}_{\max}(j\omega)$ for each system. It can be seen that Theorem \ref{thm:Small_phase_SRG} holds across the entire frequency spectrum, guaranteeing the feedback system's stability.

\subsection{Comparison between SRG-based Stability conditions} \label{ex:centralized}
We use SRGs to analyze the stability of a MIMO $4 \times 4$ system. Consider the system described by 
\begin{align*}
H_1(s)=
 \resizebox{0.8\hsize}{!}{ $ \left[\begin{array}{cccc} \dfrac{1}{\left(s+1\right)^3\,\left(s+2\right)} & \dfrac{2}{s+3} & \dfrac{4}{\left(s+1\right)\,\left(s+4\right)} & \dfrac{1}{{\left(s+1\right)}^2\,\left(s+2\right)^2}\\ 
  \dfrac{2}{s+5} & \dfrac{3}{\left(s+3\right)\,\left(s+4\right)} & \dfrac{3}{\left(s+1\right)^2\left(s+2\right)} & \dfrac{3}{s+4}\\ 
  \dfrac{1}{{\left(s+1\right)}^3} & \dfrac{3}{s+5} & \dfrac{1}{\left(s+3\right)\,\left(s+1\right)} & \dfrac{2}{\left(s+3\right)\,\left(s+4\right)}\\ 
  \dfrac{1}{s+5} & \dfrac{2}{\left(s+1\right)^5\,\left(s+2\right)} & \dfrac{1}{\left(s+1\right)\,\left(s+2\right)} & \dfrac{1}{s+1} \end{array}\right]$}
\end{align*}

\begin{figure}[h!]
    \centering
    \includegraphics[width=1\linewidth]{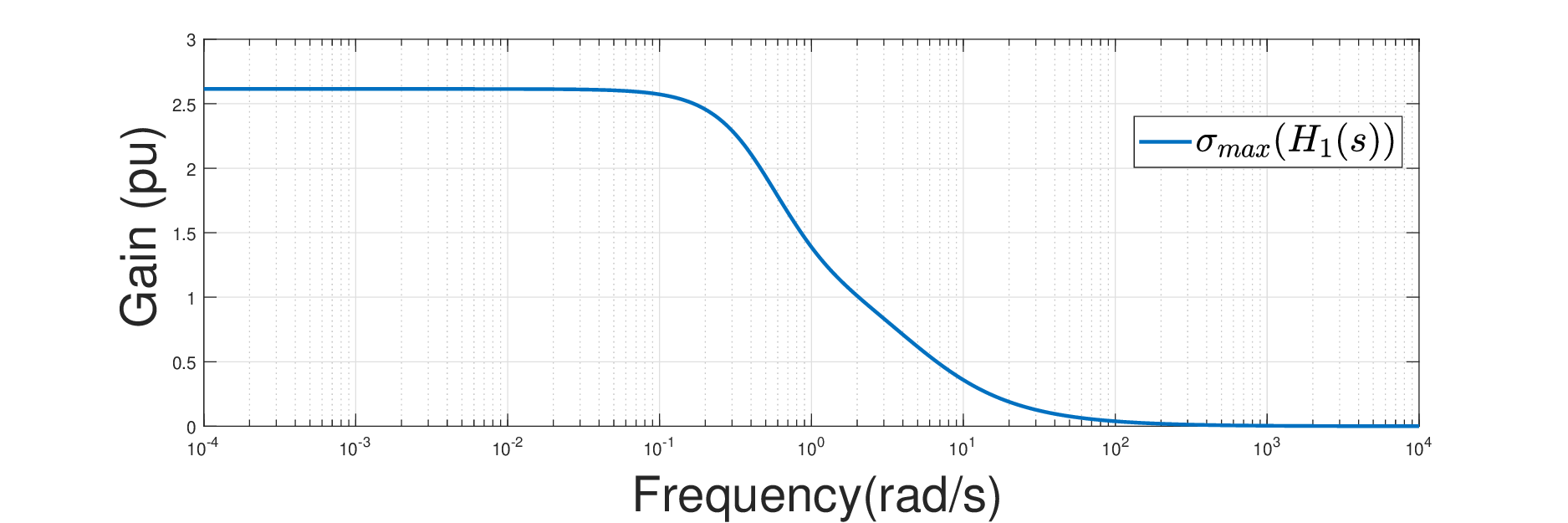}
    \caption{Gain plot $H_1(s)$}     
    \label{fig:exMIMOgain}
    \includegraphics[width=1\linewidth]{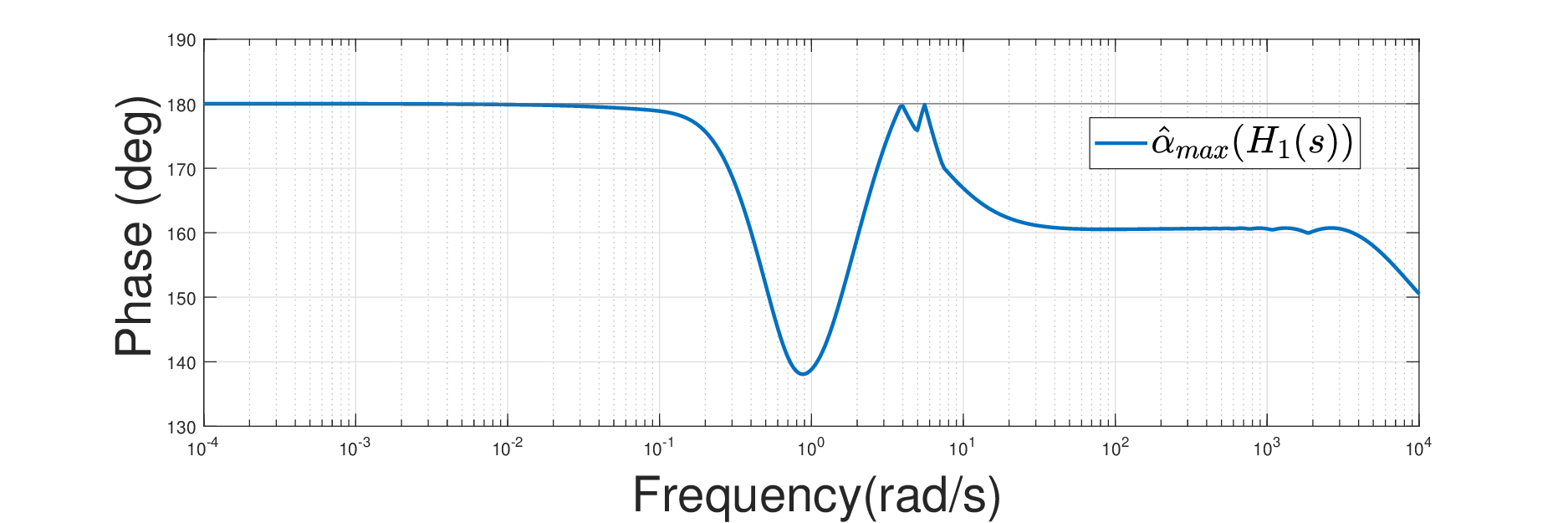}
    \caption{SRG-based phase plot $H_1(s)$}     
    \label{fig:exMIMOphase}
    \includegraphics[width=0.9\linewidth]{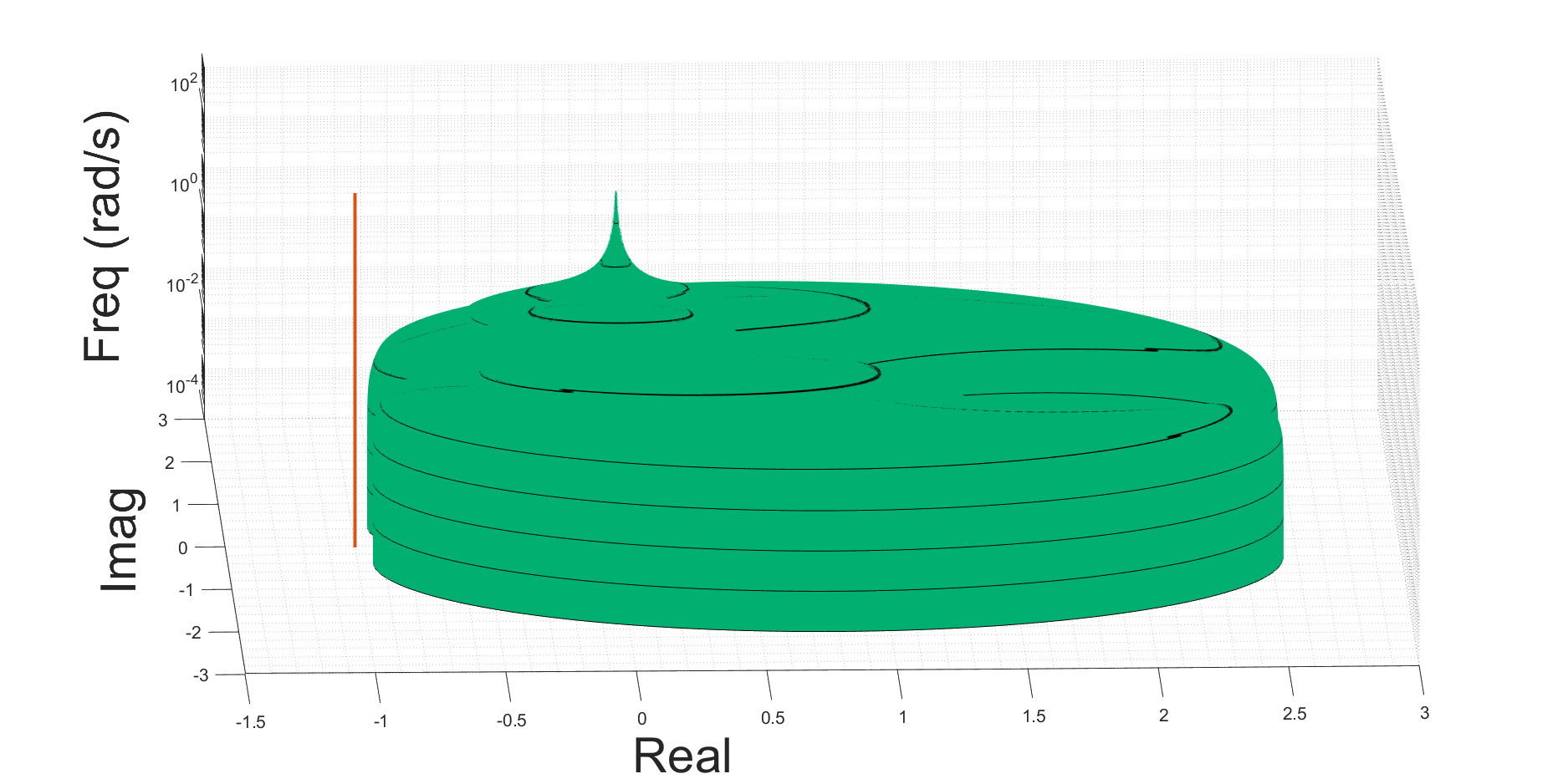}
    \caption{3D SRG$(H_1(s))$ in green with $\tau=1$ and $I_4$ in orange. Black contours are plotted to give an accurate perspective of the SRG.}         
    \label{fig:exSRG}
    \includegraphics[width=0.9\linewidth]{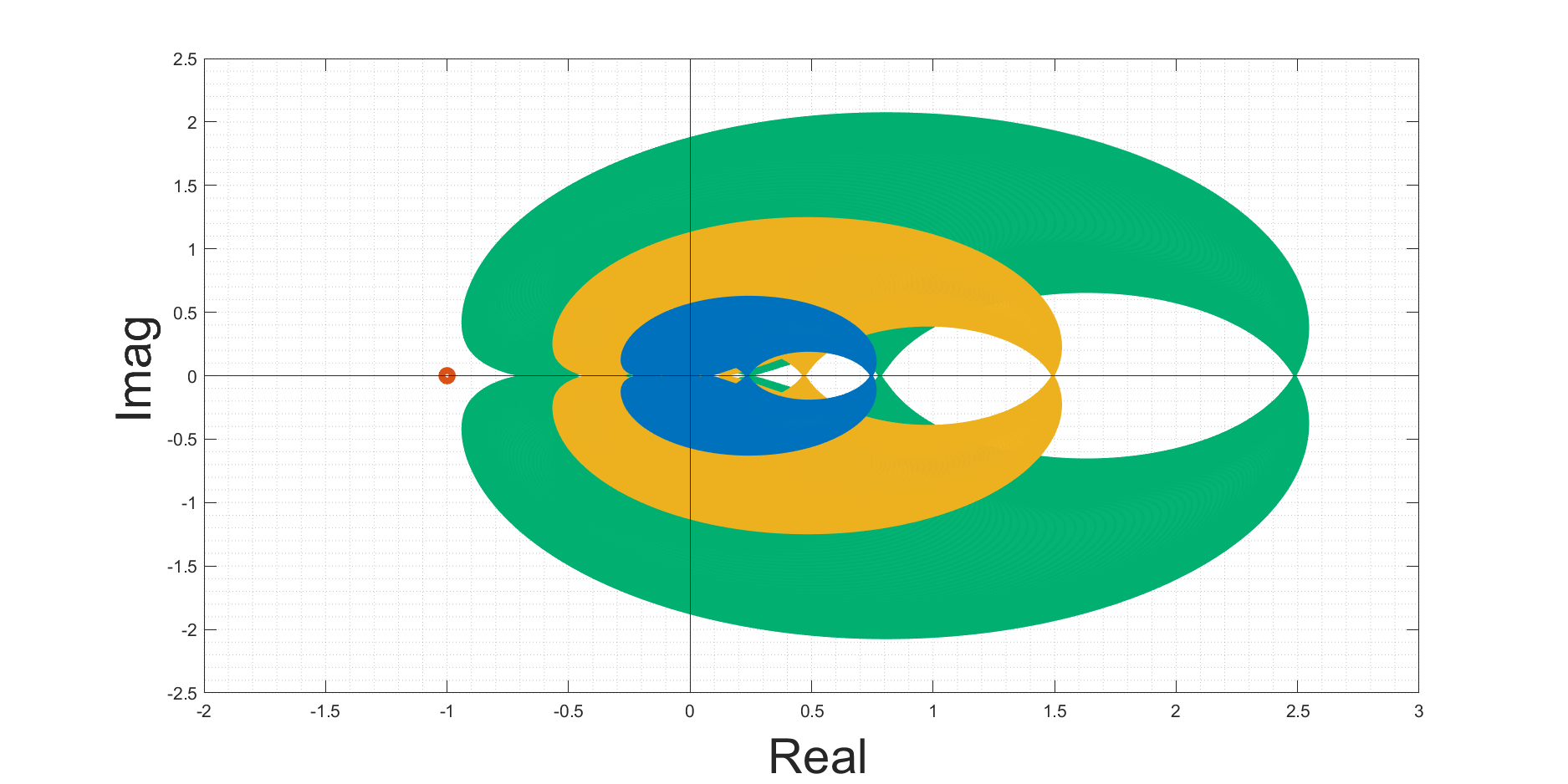}
    \caption{ 2D projection of $\SRG{H_1(s)}$ in green with $\tau=1$,  $\SRG{H_1(s)}$ in yellow with $\tau=0.6$,  $\SRG{H_1(s)}$ in blue with $\tau=0.3$. Finally, $-I_4$ in orange. }  \label{fig:exSRG_projection}
\end{figure}

Since $ H_1(s) \in \mathcal{RH}_\infty $, Corollary~\ref{corollary:nyquist} can be used to assess the stability of the negative feedback system defined by $ H_1(s) $ and $H_2(s)= I_4 $. We examine the feedback loop with the gain and phase plots shown in Figs. \ref{fig:exMIMOgain} and \ref{fig:exMIMOphase}. Figs. \ref{fig:exMIMOgain} shows that $ H_1(s) $ does not satisfy the small gain condition, particularly for frequencies below $ \omega < 0.052$~rad. Moreover, using Theorem \ref{thm:Small_phase_SRG}, the stability condition becomes $ \hat{\alpha}_{\max}(H_1(s)) < \pi $ given that $ \hat{\alpha}_{\max}(H_2(s))=0$. Nonetheless, it does not hold at all frequencies. Notably, for frequencies below $ \omega < 10^{-3} $~rad, neither the small gain nor the phase conditions are fulfilled. In this particular example, even if $H_2(s)$ is sectorial, we can not apply Theorem \ref{Thm:smallphase} because $H_1(s)$ is not sectorial at any frequency in the range $\omega\in [10^{-10},10^{10}]$~rad.

However, Corollary \ref{corollary:nyquist} can certify stability of the closed-loop system, as $ \tau\operatorname{SRG}(H(s)) \;\;\forall \tau(0,1]$ does not include the point $(-1,0)$, as depicted in Fig. \ref{fig:exSRG}. This is clearly shown in the 2D projection of the SRG in Fig.\ref{fig:exSRG_projection}, where we plot $ \tau\operatorname{SRG}(H(s))$ for $ \tau=\{1,0.6,0.3\}$. This suggests, as expected, that the approximation of the SRG by Theorem~\ref{thm:mixedpg} comes at the cost of more conservatism.
 
\section{Conclusions} \label{sec:conclusions}
We introduced two approaches for evaluating the stability of LTI feedback loop systems. The first approach leverages SRGs to develop a graphical, set-based stability criterion. In addition, we formulated a new small-phase theorem that eliminates the sectoriality assumption, a current limitation to certify system stability using this technique. As future work, we anticipate the decentralization of the small-phase theorem, further expanding its applicability.

\bibliographystyle{IEEEtran}
\bibliography{main_ECC.bib}

\begin{thebibliography}{10}
\providecommand{\url}[1]{#1}
\csname url@samestyle\endcsname
\providecommand{\newblock}{\relax}
\providecommand{\bibinfo}[2]{#2}
\providecommand{\BIBentrySTDinterwordspacing}{\spaceskip=0pt\relax}
\providecommand{\BIBentryALTinterwordstretchfactor}{4}
\providecommand{\BIBentryALTinterwordspacing}{\spaceskip=\fontdimen2\font plus
\BIBentryALTinterwordstretchfactor\fontdimen3\font minus \fontdimen4\font\relax}
\providecommand{\BIBforeignlanguage}[2]{{%
\expandafter\ifx\csname l@#1\endcsname\relax
\typeout{** WARNING: IEEEtran.bst: No hyphenation pattern has been}%
\typeout{** loaded for the language `#1'. Using the pattern for}%
\typeout{** the default language instead.}%
\else
\language=\csname l@#1\endcsname
\fi
#2}}
\providecommand{\BIBdecl}{\relax}
\BIBdecl

\bibitem{zhou1998}
K.~Zhou and J.~C. Doyle, \emph{Essentials of robust control}.\hskip 1em plus 0.5em minus 0.4em\relax Prentice hall Upper Saddle River, NJ, 1998, vol. 104.

\bibitem{1102280}
C.~Desoer and Y.-T. Wang, ``On the generalized {N}yquist stability criterion,'' \emph{IEEE Transactions on Automatic Control}, vol.~25, no.~2, pp. 187--196, 1980.

\bibitem{skogestad2005}
S.~Skogestad and I.~Postlethwaite, \emph{Multivariable feedback control: analysis and design}.\hskip 1em plus 0.5em minus 0.4em\relax john Wiley \& sons, 2005.

\bibitem{Postlethwaite1981}
I.~Postlethwaite, J.~Edmunds, and A.~MacFarlane, ``Principal gains and principal phases in the analysis of linear multivariable feedback systems,'' \emph{IEEE Transactions on Automatic Control}, 1981.

\bibitem{GRIGGS_2012_sufficientNyquist}
W.~M. Griggs, S.~S.~K. Sajja, B.~D. Anderson, and R.~N. Shorten, ``On interconnections of “mixed” systems using classical stability theory,'' \emph{Systems \& Control Letters}, 2012.

\bibitem{Chen2020}
W.~Chen, D.~Wang, and L.~Qiu, \emph{Small Phase Theorem}.\hskip 1em plus 0.5em minus 0.4em\relax London: Springer London, 2020, pp. 1--5.

\bibitem{Wang2024}
D.~Wang, W.~Chen, and L.~Qiu, ``The first five years of a phase theory for complex systems and networks,'' \emph{IEEE/CAA Journal of Automatica Sinica}, vol.~11, no.~8, pp. 1728--1743, 2024.

\bibitem{Chen2022}
W.~Chen, D.~Wang, S.~Z. Khong, and L.~Qiu, ``A phase theory of {MIMO} {LTI} systems,'' \emph{arXiv preprint arXiv:2105.03630}, 2021.

\bibitem{Chen2019}
------, ``Phase analysis of {MIMO} {LTI} systems,'' in \emph{2019 IEEE 58th Conference on Decision and Control (CDC)}, 2019, pp. 6062--6067.

\bibitem{Owens1984}
D.~Owens, ``The numerical range: a tool for robust stability studies?'' \emph{Systems \& Control Letters}, vol.~5, no.~3, pp. 153--158, 1984.

\bibitem{huang2024gain}
L.~Huang, D.~Wang, X.~Wang, H.~Xin, P.~Ju, K.~H. Johansson, and F.~Dörfler, ``Gain and phase: Decentralized stability conditions for power electronics-dominated power systems,'' 2024.

\bibitem{Woolcock2023}
L.~Woolcock and R.~Schmid, ``Mixed gain/phase robustness criterion for structured perturbations with an application to power system stability,'' \emph{IEEE Control Systems Letters}, vol.~7, pp. 3193--3198, 2023.

\bibitem{ryu2022large}
E.~Ryu and W.~Yin, \emph{Large-Scale Convex Optimization: Algorithms \& Analyses via Monotone Operators}.\hskip 1em plus 0.5em minus 0.4em\relax Cambridge Univ. Press, 2022.

\bibitem{Ryu_2021}
E.~K. Ryu, R.~Hannah, and W.~Yin, ``Scaled relative graphs: nonexpansive operators via 2d euclidean geometry,'' \emph{Mathematical Programming}, vol. 194, no. 1–2, p. 569–619, jun 2021.

\bibitem{Chaffey2021}
T.~Chaffey, F.~Forni, and R.~Sepulchre, ``Scaled relative graphs for system analysis,'' in \emph{2021 60th IEEE Conference on Decision and Control (CDC)}, 2021, pp. 3166--3172.

\bibitem{Chaffey_2023}
------, ``Graphical nonlinear system analysis,'' \emph{IEEE Transactions on Automatic Control}, vol.~68, no.~10, p. 6067–6081, Oct. 2023.

\bibitem{CHAFFEY_2022}
T.~Chaffey, ``A rolled-off passivity theorem,'' \emph{Systems \& Control Letters}, 2022.

\bibitem{arcozzi2020hardy}
N.~Arcozzi and R.~Rochberg, ``The {H}ardy space from an engineer's perspective,'' 2020.

\bibitem{psarrakos2002numerical}
P.~J. Psarrakos and M.~J. Tsatsomeros, \emph{Numerical range:(in) a matrix nutshell}.\hskip 1em plus 0.5em minus 0.4em\relax Department of Mathematics, Washington State Univ., 2002.

\bibitem{pates2021scaled}
R.~Pates, ``The scaled relative graph of a linear operator,'' \emph{arXiv preprint arXiv:2106.05650}, 2021.

\bibitem{millman1993}
R.~S. Millman and G.~D. Parker, \emph{Geometry: a metric approach with models}.\hskip 1em plus 0.5em minus 0.4em\relax Springer Science \& Business Media, 1993.

\bibitem{chaffey2024homotopy}
T.~Chaffey, A.~Kharitenko, F.~Forni, and R.~Sepulchre, ``A homotopy theorem for incremental stability,'' \emph{arXiv preprint arXiv:2412.01580}, 2024.

\bibitem{Baron2025SRG}
E.~Baron-Prada, A.~Padoan, A.~Anta, and F.~Dörfler, ``Stability results for {MIMO} {LTI} systems via scaled relative graphs,'' \emph{arXiv preprint arXiv:XXXXXXX}, 2025.

\bibitem{chen2021singularanglenonlinearsystems}
\BIBentryALTinterwordspacing
C.~Chen, W.~Chen, D.~Zhao, S.~Z. Khong, and L.~Qiu, ``The singular angle of nonlinear systems,'' 2021. [Online]. Available: \url{https://arxiv.org/abs/2109.01629}
\BIBentrySTDinterwordspacing

\bibitem{baron2023parameter}
E.~Baron-Prada, S.~Alsubaihi, K.~Alshehri, and F.~Albalawi, ``On parameter selection for first-order methods: A matrix analysis approach,'' in \emph{2023 9th International Conference on Control, Decision and Information Technologies (CoDIT)}.\hskip 1em plus 0.5em minus 0.4em\relax IEEE, 2023, pp. 445--452.

\end{thebibliography}

\appendix

\subsection{Proof Theorem \ref{thm:Small_phase_SRG}}\label{proof:Small_phase_SRG}

We start by defining the multiplication of two sets $A$ and $B$ as $    AB = \{ab | a \in A, b \in B\}$.  
Consider $H_1(j\omega)$ and  $H_2(j\omega)$ as the frequency response of the system $H_1(s)$ and  $H_2(s)$ at any arbitrary frequency $j\omega$ with $\omega \in (0,\infty]$. We recall \eqref{eqn:SRG_operator}, as 
\begin{align*}
    \operatorname{SRG}(H_1(j\omega))= \left\{\frac{\| {y}\|}{\| {x}\|} \exp \left( \pm j{\hat{\alpha}(H_1(j\omega))} \right) \right\},
\end{align*}
where $\hat{\alpha}(H_1(j\omega))=\arccos \left(\frac{\Re(\langle {y},  {x}\rangle)}{\| {y}\|\| {x}\|}\right)$. Analogously the SRG of $H_2(j\omega)$ as 
\begin{align*}
    \operatorname{SRG}(H_2(j\omega))= \left\{\frac{\| {v}\|}{\| {u}\|} \exp \left( \pm j{\hat{\alpha}(H_2(j\omega))} \right) \right\}.
\end{align*}

Then,
\begin{align*}
\operatorname{SRG}(H_1(j\omega)) \operatorname{SRG}(H_2(j\omega)) =\hspace{4cm}\\
\left\{ \dfrac{\| {y}\|}{\| {x}\|} \dfrac{\| {v}\|}{\| {u}\|} \exp \left( \pm j \hat{\alpha}(H_1(j\omega)) \pm j \hat{\alpha}(H_2(j\omega)) \right) \right\}.
\end{align*}

Since $H_1(j\omega)$ or $H_2(j\omega)$ has the right-arc property or it is over-approximated using \eqref{eqn:SRG_rightarc}, then, by \cite[Theorem 7]{Ryu_2021}.
\begin{align*}
    \operatorname{SRG}(H_1(j\omega)) \operatorname{SRG}(H_2(j\omega))=\operatorname{SRG}(H_1(j\omega) H_2(j\omega)).
\end{align*}
  In consequence, 
\begin{align}
\operatorname{SRG}(H_1(j\omega) H_2(j\omega)) =\hspace{4cm} \nonumber\\
    \resizebox{0.85\hsize}{!}{ $\left\{ \dfrac{\| {y}\|}{\| {x}\|} \dfrac{\| {v}\|}{\| {u}\|} \exp \left(  j \underbrace{(\pm \hat{\alpha}(H_1(j\omega)) \pm \hat{\alpha}(H_2(j\omega))) }_{\hat{\alpha}(H_1(j\omega) H_2(j\omega))} \right) \right\}$}. 
    \label{eqn:symmetryangles}
\end{align}
Using \eqref{eqn:SingularAngledef}, the maximum phase of  $\operatorname{SRG}(H_1(j\omega) H_2(j\omega))$ can be rewritten as
\begin{align}
      \hat{\alpha}_{\max}(H_1(j\omega) H_2(j\omega))=  \hat{\alpha}_{\max}(H_1(j\omega)) + \hat{\alpha}_{\max}(H_2(j\omega)).
      \label{eqn:smallphase_step1}
\end{align} 
Thus, using \eqref{eqn:SRGsmallphase}  to bound \eqref{eqn:smallphase_step1} we obtain
\begin{align}
     \hat{\alpha}_{\max}(H_1(j\omega)& H_2(j\omega))= \nonumber\\
     &\hat{\alpha}_{\max}(H_1(j\omega)) + \hat{\alpha}_{\max}(H_2(j\omega)) < \pi \label{eqn:Smallphasefinal}.
\end{align}


Equation \eqref{eqn:Smallphasefinal} states that the maximum angle of $\operatorname{SRG}(H_1(j\omega) H_2(j\omega))$, $\hat{\alpha}_{\max}(H_1(j\omega) H_2(j\omega))$, must fall strictly within the $(- \pi, \pi)$ range. This provides a sufficient condition to satisfy Corollary \ref{corollary:nyquist}, which ensures that $(-1, 0)$ does not lie within $\tau \operatorname{SRG}(H_1(j\omega) H_2(j\omega))$ for any $\tau \in (0, 1]$, thus guaranteeing feedback system stability.
$\hfill\blacksquare$

\subsection{ Proof Theorem \ref{thm:mixedpg}}\label{proof:mixedpg}

We can assess the stability of the closed-loop system on a frequency-by-frequency basis by applying the superposition principle \cite{zhou1998}, for each $ \omega \in [0, \infty) $. Exploiting this principle, we start by defining the following two sets
\begin{align*}
    \omega_{\hat{\alpha}}&=\{\omega\; |\omega \in [0,\infty), 
    \hat{\alpha}_{\max}(H_1(j\omega)) + \hat{\alpha}_{\max}(H_2(j\omega)) < \pi\}, \\
    \omega_\sigma&=\{\omega\; |\omega \in [0,\infty), 
     \sigma_{\max}(H_1(j\omega))\sigma_{\max}(H_2(j\omega))<1 \},
\end{align*}
 where $\omega_\alpha$ and $\omega_\sigma$ are the set of frequencies that meet Theorem~\ref{thm:Small_phase_SRG} and Theorem~\ref{theorem:smallgain}, respectively. If the union of both sets covers the entire frequency spectrum, i.e.,
 \begin{align}
     \omega_{\hat{\alpha}} \cup \omega_\sigma = [0,\infty),
 \end{align}
 Then $\operatorname{SRG}(H_1(j\omega) H_2(j\omega))$ does not include $(-1, 0)$ and thus the feedback system is stable.
 $\hfill\blacksquare$

\end{document}